\documentclass[12pt,oneside]{article}
\usepackage{amsthm,verbatim,mathtools}
\usepackage{amssymb,latexsym}
\usepackage{amsmath}
\usepackage{graphicx}
\usepackage{color,cite}
\usepackage[dvipsnames]{xcolor}
\usepackage{cancel}
\usepackage{soul}
\usepackage{enumerate}
\newtheorem{theorem}{Theorem}[section]

\newtheorem{conjecture}[theorem]{Conjecture}
\newtheorem{proposition}[theorem]{Proposition}

\begin{document}
\title{On a conjecture about a class of permutation trinomials}
\author{Daniele Bartoli}
\date{}
\maketitle

\begin{abstract}
We prove a conjecture by Tu, Zeng,  Li, and Helleseth concerning trinomials $f_{\alpha,\beta}(x)= x + \alpha x^{q(q-1)+1} + \beta x^{2(q-1)+1} \in \mathbb{F}_{q^2}[x]$, $\alpha\beta \neq 0$, $q$ even, characterizing all the pairs $(\alpha,\beta)\in \mathbb{F}_{q^2}^2$ for which $f_{\alpha,\beta}(x)$ is a permutation of $\mathbb{F}_{q^2}$. 
\end{abstract}

 
\section{Introduction}

Let $q=p^h$ be a prime power. A polynomial $f(x)\in \mathbb{F}_q[x]$ is a {\it permutation polynomial} (PP for short) if it is a bijection of the finite field $\mathbb{F}_q$ into itself. On the other hand, each permutation of $\mathbb{F}_q$ can be expressed as a polynomial over $\mathbb{F}_q$. 

Permutation polynomials have nice connections with applied areas of mathematics, such as cryptography, coding theory, and combinatorial designs. Random PP for a given field $\mathbb{F}_q$ can be easily constructed. In many  applications, however, simple structures or additional extraordinary properties on PPs are usually required and PPs meeting these criteria are usually difficult to find. For a deeper treatment of the connections of PPs with other fields of mathematics we refer to \cite{MuPa, Hou2015} and the references therein.

In this work we deal with a particular class of PP{s} introduced in  \cite{TZLH2017}, that is polynomials $f_{\alpha,\beta}(x)\in \mathbb{F}_{q^2}$ of  type
\begin{equation}\label{f_alfa_beta}
 x + \alpha x^{q(q-1)+1} + \beta x^{2(q-1)+1},
\end{equation}
with $\alpha ,\beta\in \mathbb{F}_{q^2}^*$, $q=2^m$. The authors prove that if 
\begin{enumerate}
\item $\beta=\alpha^{q-1}$ and $Tr\left(1+\frac{1}{\alpha^{q+1}}\right)=0$ or
\item $\beta(1+\alpha^{q+1}+\beta^{q+1})+\alpha^{2q}=0$, $\beta^{q+1}\neq 1$, and $Tr\left( \frac{\beta^{q+1}}{\alpha^{q+1}}\right)=0$
\end{enumerate}
then $f_{\alpha,\beta}(x)$ permutes $\mathbb{F}_{q^2}$; see \cite[Theorem 1]{TZLH2017}. Due to some experimental results, the authors state the following conjecture.

\begin{conjecture}\label{Conjecture}
If $f_{\alpha,\beta}(x)$ permutes $\mathbb{F}_{q^2}$ then $1.$ or $2.$ holds.
\end{conjecture}

In this work we prove the above conjecture, using the well known connection between permutation polynomials and  algebraic curves over finite fields. 

First of all, let us remark that  the polynomials $f_{\alpha,\beta}(x)$ belong to the more general class of polynomials 
$$f_{r,d,h}(x)=x^rh\left( x^{\frac{q-1}{d}}\right),$$
where  $h(x) $ is a polynomial over  $\mathbb{F}_{q}$, $q=p^m$, $d$ a divisor of $q-1$,  $r$ an integer with  $1 \leq r<(q-1)/d$.

A useful criterion to decide weather $f_{r,d,h}$ permutes $\mathbb{F}_{q}$ was established in \cite{PL2001,Zieve2009}.

\begin{theorem}\label{ThPLZ} 
The polynomial $f_{r,d,h}(x)$ is a PP of $\mathbb{F}_{q}$ if and only if $gcd(r, (q-1)/d) =1$ and $x^rh(x)^{(q-1)/d}$ permutes the set $\mu_d$ of the $d$-th roots of unity in $\mathbb{F}_{q}$.
\end{theorem}

The above theorem can be seen as an application of the AGW criterion; see \cite{AGW2011,YD2011,YD2014}.

By Theorem \ref{ThPLZ}, the polynomial $f_{\alpha,\beta}(x)$ in \eqref{f_alfa_beta} permutes $\mathbb{F}_{q^2}$ if and only if  
$$x(1 + \alpha x^{q} + \beta x^{2})^{q-1}$$
permutes $\mu_{q+1}$. Recalling that for $x\in \mu_{q+1}$ we have that $x^q=1/x$, this is equivalent to 
$$g_{\alpha,\beta}(x)=\frac{\alpha^q x^3+x^2+\beta^q}{\beta x^3+x+\alpha}$$
permuting $\mu_{q+1}$. 

We call the rational function $g_{\alpha,\beta}(x)$ the \emph{fractional polynomial} associated with the polynomial $f_{\alpha,\beta}(x)$. Permutational properties of fractional polynomial have been studied in several works. Algebraic curves associated with fractional polynomials have been investigated in \cite{BG2017}. 

Consider the fractional polynomial $g_{\alpha,\beta}(x)$ and let $z\in \mathbb{F}_{q^2}\setminus \mathbb{F}_{q}$ be fixed. Every element in $\mu_{q+1}\setminus \{1\}$ can be written as $\varphi(x)=\frac{x+z}{x+z^q}$ for some $x \in \mathbb{F}_q$. Define $\infty$ to be an element such that $\varphi(\infty)=1$. Therefore the function $G(x)=g_{\alpha,\beta}(\varphi(x)) : \mathbb{F}_q\cup \{\infty\} \to \mu_{q+1}$ is a bijection if and only if $g_{\alpha,\beta}(x)$ is a bijection of $\mu_{q+1}$. This happens if and only if  $F(x)=G(x)_{\mid\mathbb{F}_q} : \mathbb{F}_q \to \mu_{q+1}\setminus\{(\alpha+\beta+1)^{q-1}\}$ is a bijection. 

Let $\mathcal {C}_{F}$ be the curve of affine equation $(F(x)-F(y))/(x-y)=0$. It is easily verified that $\mathcal{C}_{F}$ is defined over $\mathbb{F}_q$, and therefore $F(x)$ is a bijection if and only if $\mathcal{C}_{F}$ has no $\mathbb{F}_q$-rational points off the line $x=y$. 

In general, to establish if a particular curve $\mathcal{C}$ does not have $\mathbb{F}_{q}$-rational points is a  hard problem.  If the degree $d$ of  $\mathcal{C}$ is small with the respect to the size $q$ of the field $\mathbb{F}_q$, namely $d<\sqrt[4]{q}$, then the existence of an absolutely irreducible $\mathbb{F}_q$-rational component in $\mathcal{C}$ distinct from $x=y$ yields, due to Hasse-Weil Theorem (see \cite[Theorem 5.2.3]{Sti})), the existence of $\mathbb{F}_q$-rational points off the line $x=y$. 

In our case (see Section \ref{Sec:Curva}) the curve $\mathcal{C}_F$ has degree at most $4$. For $q$ large enough such a curve has no $\mathbb{F}_q$-rational points off $x=y$ if only if it splits into absolutely irreducible components not defined over $\mathbb{F}_q$ which have no $\mathbb{F}_q$-rational points off the line $x=y$. This gives us a useful criterion to decide whether the fractional polynomial $g_{\alpha,\beta}$ permutes $\mathbb{F}_{q^2}$.

We would like to point out that our results are valid for general $q$'s of type $q=2^m$, and obviously do not rely on computer searches carried out in $\mathbb F_q$ for large $q$. However, our method requires a computer to assist us in computing resultants between polynomials and in factorizing  polynomials of low degrees over small fields.  The elementary MAGMA \cite{MAGMA} programs used for our purposes are presented in Appendix.

The paper is organized as follows. In Section \ref{Sec:Curva} we study the curve $\mathcal{C}_{F}$ associated with the fractional polynomial  $g_{\alpha,\beta}$. In particular we determine  necessary and sufficient conditions for the corresponding curve to split completely into absolutely irreducible components not defined over $\mathbb{F}_{q}$. In Section \ref{Sec:Proof} we prove Conjecture \ref{Conjecture}. Finally, Appendix contains the short MAGMA \cite{MAGMA} programs we used in Section \ref{Sec:Curva}.

\section{Investigation of the curve $\mathcal{C}_{F}$}\label{Sec:Curva}
In this section we investigate the curve $\mathcal{C}_{F}$ associated with the fractional polynomial  $g_{\alpha,\beta}$. Our goal is to determine all the possible pairs $(\alpha,\beta) \in \mathbb{F}_{q^2}$ for which the curve $\mathcal{C}_{F}$ splits completely into absolutely irreducible components not defined over $\mathbb{F}_q$. 

Let $q=2^m$. Denote by $Tr: \mathbb{F}_q \to\mathbb{F}_2$ the absolute trace defined by  
$$Tr(z)=z^{2^{m-1}}+z^{2^{m-2}}+\cdots + z^{2^2}+z^{2}+z$$
for any $z\in \mathbb{F}_{q}$. Let $k\in \mathbb{F}_q$ have trace equal to $1$. Consider $i\in \mathbb{F}_{q^2}$  such that $i^2=i+k$. In particular $i^q=i+1$. As mentioned above, since $i \in \mathbb{F}_{q^2}\setminus \mathbb{F}_q$,  each $\overline{x} \in \mu_{q+1} \setminus\{1\}$ can be written as $\frac{x+i}{x+i+1}$ for some $x\in \mathbb{F}_q$. Also, $\{1,i\}$ is a basis of $\mathbb{F}_{q^2}$ over $\mathbb{F}_q$ and every element $z$ in $\mathbb{F}_{q^2}$ can be written as $z=z_1+iz_2$ for some $z_1,z_2 \in \mathbb{F}_q$. 

Let $\alpha=A+iB$ and $\beta=C+iD$, $A,B,C,D\in \mathbb{F}_q$. We already mentioned that the  fractional polynomial $g_{\alpha,\beta}$ is a permutation of $\mu_{q+1}$ if and only if  the curve $\mathcal{C}_F$ 
has no $\mathbb{F}_q$-rational points off the line $x=y$. 
By direct computations (see also Appendix), such a curve has affine equation $L(x,y)=0$, where
\begin{equation}\label{EquazioneCurva}
L(x,y)=\gamma_{2,2}x^2y^2+ \gamma_{2,1}x^2y+\gamma_{1,2}xy^2+\gamma_{2,0}x^2+\gamma_{1,1}xy+\gamma_{0,2}y^2+\gamma_{1,0} x+\gamma_{0,1}y+\gamma_{0,0},
\end{equation}
with 
$$\begin{array}{lll}
\gamma_{2,2}&=&A^2 + A B + B^2 k + C^2 + C D + D^2 k + D + 1,\\
\gamma_{2,1}&=&A^2 + A B + B^2 k + B + C^2 + C D + D^2 k + 1,\\
\gamma_{1,2}&=&A^2 + A B + B^2 k + B + C^2 + C D + D^2 k + 1,\\
\gamma_{2,0}&=&A^2 k + A^2 + A B k + A B + A + B^2 k^2 + B^2 k + C^2 k \\
		     &&+ C^2 + C D k + C D + C + D^2 k^2 + D^2 k + D k + D + k,\\
\gamma_{0,2}&=&A^2 k + A^2 + A B k + A B + A + B^2 k^2 + B^2 k + C^2 k \\
		     &&+ C^2 + C D k + C D + C + D^2 k^2 + D^2 k + D k + D + k,\\
\gamma_{1,1}&=&A^2 + A B + B^2 k + C^2 + C D + D^2 k + 1,\\
\gamma_{1,0}&=&A^2 k + A B k + A + B^2 k^2 + B k + C^2 k + C D k + D^2 k^2 + k,\\
\gamma_{0,1}&=&A^2 k + A B k + A + B^2 k^2 + B k + C^2 k + C D k + D^2 k^2 + k,\\
\gamma_{0,0}&=&A^2 k^2 + A B k^2 + B^2 k^3 + C^2 k^2 + C D k^2 + C + D^2 k^3 + D k^2 + D k + D + k^2.\\
\end{array}
$$

In what follows  we determine all the  cases in which $\mathcal{C}_{F}$ splits into components none of them is defined over $\mathbb{F}_q$. 

\begin{proposition}\label{Proposizione:riassunto}
Let $k$ be a fixed element in $\mathbb{F}_q$ with $Tr(k)=1$ and $i\in \mathbb{F}_{q^2}\setminus \mathbb{F}_q$ be such that $i^2=i+k$. Also, let $\overline{k}$ be such that $\overline{k}^4=k$ and consider $\alpha=A+iB$, $\beta=C+iD$, with $A,B,C,D \in \mathbb{F}_q$ and $\alpha\beta\neq 0$. The curve $\mathcal{C}_F$ splits into absolutely irreducible components not defined over $\mathbb{F}_q$ if and only if one of the following cases occurs.  
\begin{enumerate}
\item $$A=\xi^2+\xi, \qquad B=D=0, \qquad C=\xi^2,$$
where $\xi\in \mathbb{F}_q\setminus\{0,1\}$, $Tr\left(\xi/(\xi+1)\right)=0$.

\item $$A=(\eta \overline{k} + \eta + \xi)^2(\eta^2\overline{k} + \eta \xi + \xi^2 + 1), \quad B=  \eta^4\overline{k}+ \xi^2 \eta^2+ \xi \eta^3+ \eta^2, \quad C=\xi^4, \quad D=\eta^4,$$
 where $\xi,\eta \in \mathbb{F}_q$, $\eta\ne 0$, $Tr\left(\frac{B}{D}  + 1 +\frac{1}{D^2}  +  \frac{D}{B^2}\right)=1$, $\eta^2\overline{k}+ \xi^2+ \xi \eta+ 1\neq 0, \eta^2$.

\item $$A\in \mathbb{F}_{q}^*, \qquad B=D=0, \qquad C=1,$$
with $Tr \left(1+\frac{1}{A}\right)= 0$.

\item $$A=\frac{BC+BD+B}{D},\qquad B,C \in \mathbb{F}_q, \qquad D \in \mathbb{F}_q^{*},$$
with $BD(B+D)\neq 0$, $kD^2+C^2+CD+1=0$, $Tr\left(1+\frac{D}{B^2}\right)=0$.
 
\item  $q=2^{2s+1}$, 
$$C=A+B+1, \qquad D=B,$$
with  $A^2 + A B + B^2 k + B = 0$ and $B \neq 0$.

\end{enumerate}

\end{proposition}

The proof of this result is given in the following subsections. We distinguish the cases  $\gamma_{2,2}\neq 0 $ and $\gamma_{2,2}=0$.

\subsection{Case $\gamma_{2,2}\neq 0$}
In this case $\mathcal{C}_F$ has degree $4$. We observe that the morphism $(x,y)\to (y,x)$ fixes $\mathcal{C}_F$ and therefore it acts on its components. Also, note that if $\mathcal{C}_F$ spits completely into absolute irreducible components not defined over $\mathbb{F}_q$ then the only possibilities are 2 conics or 4 lines. 

\begin{enumerate}
\item $\mathcal{C}_F$ splits into four lines. In this case the factorization of $L(x,y)$ in \eqref{EquazioneCurva} must be 
$$\gamma_{2,2}(x+a)(x+b)(y+a)(y+b)$$
 for some $a,b \in \overline{\mathbb{F}_q}$, since the homogeneous part of degree $4$ in $\mathcal{C}_F$ is $\gamma_{2,2}x^2y^2$.

{\bf Suppose $B=D$.} Since $A^2 + A B + B^2 k + C^2 + C D + D^2 k + D + 1\neq 0$, we get $(A + C + 1)(A + B + C + 1)\neq 0$. 
Comparing the coefficients, we obtain in particular $B=0$ and $A^2 + C^2 + C=0$. Since $A\neq 0$ otherwise $\alpha=A+iB=0$,  we have that $C=\xi^2$, $A=\xi^2+\xi$ for some $\xi \in \mathbb{F}_{q}\setminus\{0,1\}$. On the other hand, if $B=D=0$, $C=\xi^2$, $A=\xi^2+\xi$, for some $\xi \in \mathbb{F}_q \setminus \{0,1\}$ then $\mathcal{C}_F$ splits over $\mathbb{F}_q$ in  
$$(\xi x^2 + x^2 + \xi x + x + \xi k + \xi  + k)(\xi y^2 + y^2 + \xi y + y + \xi k + \xi + k)=0.$$ 
The four lines are not defined over $\mathbb{F}_{q}$ if and only if the roots of the polynomial 
$$\xi T^2 + T^2 + \xi T + T + \xi k + \xi  + k$$
are not in $\mathbb{F}_q$. This happens if and only if (see for instance \cite[Section 1.4]{HirschBook})
$$Tr\left( \frac{(\xi+1)(\xi k + \xi  + k)}{(\xi+1)^2} \right)=1 \iff Tr\left( \frac{\xi k + \xi  + k}{\xi+1} \right)=1$$
$$\iff Tr\left( k+\frac{\xi}{\xi+1}\right)=1 \iff Tr\left(\frac{\xi}{\xi+1}\right)=0 \iff Tr\left(\frac{C}{C+1}\right)=0.$$

{\bf Suppose now $B\neq D$.} Then in particular 

$$A^2D + ABD + B^2Dk + B^2 + C^2D + CD^2 + D^3k + D=0.$$
So, $D\neq 0$ otherwise $B=0=D$. Finally, from  $A D^2 + B^3 + B C D + B D^2 + B D=0$
we have
\begin{equation}\label{Cond1}
A=\frac{B^3 + B C D + B D^2 + B D}{D^2}, \qquad k =\frac{   B^4 + C^2D^2 + CD^3  + D^2 }{D^4}.
\end{equation}
Let $ C=\xi^4$, and  $D=\eta^4$, then 
 $$A=(\eta \overline{k} + \eta + \xi)^2(\eta^2\overline{k} + \eta \xi + \xi^2 + 1), \qquad B=  \eta^2(\eta^2\overline{k}+ \xi^2+ \xi \eta+ 1).$$
 Note that $B\neq 0$ otherwise, from $D\neq 0$, then we get $A=0$, a contradiction. 
On the other hand, if \eqref{Cond1} holds then the polynomial $L(x,y)$ splits in 
$$((D^5+BD^4)y^2 + BD^4y  + B^5 + B^4 D + B^2 D^3 + B C^2D^2 + B C D^3 + B D^2 + C^2 D^3 + D^5 + D^4 + D^3)\cdot$$
$$((D^5+BD^4)x^2 + BD^4x  + B^5 + B^4 D + B^2 D^3 + B C^2D^2 + B C D^3 + B D^2 + C^2 D^3 + D^5 + D^4 + D^3)=0.$$
The four lines are not defined over $\mathbb{F}_{q}$ if and only if 

${\small Tr\left(\frac{B^6 D^4 + B^4 D^6 + B^3 D^7 + B^2 C^2 D^6 + B^2 C D^7 + B^2 D^8+B^2 D^6 + B C D^8 + B D^9 + B D^8 + C^2 D^8 + D^{10} + D^9 + D^8}{B^2D^8}\right)=1}$\\
$\iff {\small Tr\left(\frac{B^4}{D^4} + \frac{B^2}{D^2} + \frac{B}{D} + \frac{C^2}{D^2} + \frac{C}{D} + 1 +\frac{1}{D^2} + \frac{C}{B} + \frac{D}{B} + \frac{1}{B} + \frac{C^2}{B^2} + \frac{D^2}{B^2} + \frac{D}{B^2} + \frac{1}{B^2}\right)=1}$\\
$$\iff Tr\left(\frac{B}{D}  + 1 +\frac{1}{D^2}  +  \frac{D}{B^2}\right)=1.$$

\item $\mathcal{C}_F$ splits into two absolute irreducible conics not defined over $\mathbb{F}_q$ both fixed by $(x,y)\mapsto (y,x)$. The two conics must be switched by $(x,y) \mapsto (x^q,y^q)$ and therefore  the factorization of $L(x,y)$ must be 
$$(x y+(\sigma+i \tau) x+(\sigma+i \tau) y+(\gamma+i \delta))(x y+(\sigma+\tau+i \tau) x+(\sigma+\tau+i \tau) y+(\gamma+\delta +i \delta)),$$
  for some $\sigma,\tau,\gamma,\delta \in \mathbb{F}_q$.  In particular (see Appendix) we have that 

$$
\begin{array}{lll}
z_2 \sigma^2+z_1\sigma+z_0=0\\
\\
\delta =\frac{ A^2 + A B + B^2 k + C^2 + C D + D^2 k + 1}{A^2 + A B + B^2 k + C^2 + C D + D^2 k + D + 1}\\
\\
\tau =\frac{A^2 + A B + B^2 k + B + C^2 + C D + D^2 k + 1}{A^2 + A B + B^2 k + C^2 + C D + D^2 k + D + 1}\\
\\
\gamma =\frac{(A^2 + A B + B^2 k + C^2 + C D + D^2 k + 1)a+A^2 k + A B k + A + B^2 k^2 + B k + C^2 k + C D k + D^2 k^2 + k}{A^2 + A B + B^2 k + C^2 + C D + D^2 k + D + 1},\\
\\
H(A,B,C,D)=0,\\
\end{array}
$$
where
$$\begin{array}{lll}  
z_0&=&A^4 + A^3 + A^2 B^2 + A^2 B + A^2 C + A^2 + A B^2 k + A B C + \\
&&A B + A C^2 +  A C D + A D^2 k + A D + A + B^4 k^2 + B^2 C k + C^4 + C^3\\
&& + C^2 D^2 +  C^2 D + C^2 + C D^2 k + C + D^4 k^2 + D^2 + D\\
 z_1&=&A^4 + A^2 B^2 + A^2 B + A^2 D + A B^2 + A B D + B^4 k^2 + B^3 k \\
&&+ B^2 D k +  B C^2 + B C D + B D^2 k + B D + B + C^4 + C^2 D^2 + C^2 D\\
&& + C D^2 +  D^4 k^2 + D^3 k + D + 1\\
z_2&=& A^4 + A^2 B^2 + B^4 k^2 + C^4 + C^2 D^2 + D^4 k^2 + D^2 + 1,\\
\end{array}
$$    
and 
$$
\begin{array}{ll}
H(A,B,C,D) =&A^6 + A^5 B + A^4 B^2 k + A^4 B^2 + A^4 C^2 + A^4 C D + A^4 D^2 k \\
		    &+ A^4 + A^3 B^3 + A^2 B^4 k^2 + A^2 B^4 k + A^2 B^2 C^2 + A^2 B^2 C D\\
		    & + A^2 B^2 D^2 k + A^2 B^2 +  A^2 C^4 + A^2 C^2 D^2 + A^2 C^2 + A^2 C D\\
		    & + A^2 D^4 k^2 + A^2 D^2 k + A^2 D +  A B^5 k^2 + A B C^4 + A B C^2 D^2\\
		    & + A B C^2 + A B C D + A B D^4 k^2 + A B D^2 k +  B^6 k^3 + B^4 C^2 k^2\\
		    & + B^4 C D k^2 + B^4 D^2 k^3 + B^4 k^2 + B^2 C^4 k + B^2 C^2 D^2 k\\
		    & + B^2 C^2 k + B^2 C D k + B^2 C + B^2 D^4 k^3 + B^2 D^2 k^2 + B^2 D k\\
		    & + B^2 D + C^6 + C^5 D + C^4 D^2 k + C^4 D^2 + C^3 D^3 + C^2 D^4 k^2\\
		    & + C^2 D^4 k + C^2 +C D^5 k^2 + C D + D^6 k^3 + D^2 k.\\
\end{array}
$$
By direct checking (see Appendix) if $H(A,B,C,D)=0$ then 
$$\Big((\sigma+i \tau)^2+(\gamma+i \delta)\Big) \cdot \Big((\sigma+\tau +i \tau)^2+(\gamma+\delta +i \delta)\Big)=0,$$
that is at least one of the two conics splits into two lines, namely 
$$\Big(x +(\sigma+i \tau)\Big)\cdot \Big(y +(\sigma+i \tau)\Big)$$
or 
$$\Big(x +(\sigma+\tau+i \tau)\Big)\cdot \Big(y +(\sigma+\tau+i \tau)\Big).$$
This is not possible since by assumption the two conics are absolutely irreducible.

\item $\mathcal{C}_F$ splits into two absolute irreducible conics not defined over $\mathbb{F}_q$ switched  by $(x,y)\mapsto (y,x)$. The two conics must be switched also by $(x,y) \mapsto (x^q,y^q)$ and therefore  the factorization of $L(x,y)$ must be either 
$$(x^2+(\sigma+i \tau)x+(\sigma+\tau +i \tau) y+c) (y^2+(\sigma+\tau +i \tau) x+(\sigma+i \tau) y+c)$$ 
or 
$$(x y+(\sigma+i \tau) x+(\sigma+\tau +i \tau) y+c)(x y+(\sigma+\tau +i \tau) x+(\sigma+i \tau) y+c),$$  for some $\sigma, \tau,c \in \mathbb{F}_q$. In the former case, comparing the coefficients, we immediately get $\sigma=\tau=0$, that is $\mathcal{C}_F$ splits into four line. 

Consider now the latter case.
\begin{itemize} 
\item $B=D$. 
As above $(A + C + 1)(A + B + C + 1)\neq 0$. This implies $\tau=1$ and $B=D=0$. Also $A^2 k + A^2 c + A + C^2 k + C^2 c + k + c=0$ and $(C+1)(A^2 + C^2 + C)=0$. If $A^2 + C^2 + C=0$ then the curve splits into four lines (see above). So $C=1$, $A \sigma^2 + A \sigma + A + 1=0$, and $A k + A c + 1 =0$. Therefore 
$$Tr\left(1+\frac{1}{A}\right)=0.$$

\item $B\neq D$. Then in particular  
$$\begin{array}{lll}
(A^2 + A B + B^2 k + C^2 + C D + D^2 k + D + 1)\tau\\
\hspace{ 2 cm}+A^2 + A B + B^2 k + B + C^2 + C D + D^2 k + 1&=&0\\
\\
(B^2 D + D^3)k+A^2 D + A B D + B^2 + C^2 D + C D^2 + D&=&0\\
\\
(B^3 D + B D^3)c+A^2 B D + A B^2 D + A B D^2\\
\hspace{2 cm} + A D^3 + B^3 + B C^2 D + B C D^2 + B D&=&0\\
\\
D(A D + B C + B D + B)(A D^2 + B^3 + B C D + B D^2 + B D)&=&0.\\
\end{array}
$$
It is easily seen that $D=0$ implies $B=0$, a contradiction. Also, if $B=0$, since $A$ must be different from $0$, then $D=0$, again a contradiction. 

Therefore either 
$$A D + B C + B D + B=0\qquad  \textrm{ or } \qquad A D^2 + B^3 + B C D + B D^2 + B D=0.$$

If $A=\frac{B C + B D + B}{D}$, then $k=\frac{C^2 + C D  + 1}{D^2}$, $\tau=\frac{B}{B+D}$, $c=\frac{BC^2 + BCD +   B + C^2D   + D^3 + D^2 + D}{BD^2+D^3}$ and 
$$(B^2 + D^2)\sigma^2+
(B^2 + BD)\sigma+
(B^2 + BC + BD + B + C^2 + D^2 + D + 1)=0.$$
Since $\sigma\in \mathbb{F}_q$, we have 
 
$$Tr\left(\frac{(B^2 + D^2)(B^2 + BC + BD + B + C^2 + D^2 + D + 1)}{(B^2 + BD)^2}\right)=Tr\left(1+\frac{D}{B^2}\right)=0.$$

If $A=\frac{B^3 + B C D + B D^2 + B D}{D^2}$ then $k=\frac{B^4 + C^2 D^2 + C D^3  + D^2}{D^4}$ and the polynomial $L(x,y)$ factorizes in 
$$((B D^4 + D^5)y^2 + B D^4y + B^5 + B^4 D + B^2 D^3 + B C^2 D^2 + B C D^3
        + B D^2 + C^2 D^3 + D^5 + D^4 + D^3)\cdot$$
$$    ((B D^4 + D^5)x^2 + B D^4x + B^5 + B^4 D + B^2 D^3 + B C^2 D^2 + B C D^3
        + B D^2 + C^2 D^3 + D^5 + D^4 + D^3)$$
that is the curve contains four lines, a contradiction.
\end{itemize}
\end{enumerate}

\subsection{Case $\gamma_{2,2}= 0$}

In this case the curve $\mathcal{C}_F$ has degree less than $4$. Again we want to determine when $\mathcal{C}_F$ contains no absolute irreducible components defined over $\mathbb{F}_q$. The only possible case is when $\mathcal{C}_F$ factorizes into lines. 

Suppose $B\neq D$. In this case the curve $\mathcal{C}_F$ has degree exactly three since the coefficient of both $x^2y$ and $xy^2$ is 
$$A^2 + A B + B^2 k + B + C^2 + C D + D^2 k + 1,$$
which is different from zero. The homogeneous part of the polynomial $L(x,y)$ is then $xy(x+y)$. Since the ideal point of the line $x=0$ is a simple point for $\mathcal{C}_F$ and it is clearly $\mathbb{F}_q$-rational, there exists exactly one absolutely irreducible component of $\mathcal{C}_F$ through it. Then, such a component must be fixed by the morphism $(x,y) \mapsto (x^q,y^q)$ and therefore it is $\mathbb{F}_q$-rational. Since it cannot be the line $x+y=0$ the corresponding function $F$ is not a bijection.

So we must suppose $B=D$ and then from $A^2 + A B + B^2 k + C^2 + C D + D^2 k + D + 1=0$ we get $A^2 + A B  + C^2 + C B  + B + 1=(A + C + 1)(A + B + C + 1)=0$. 

If $A + C + 1=0$ then the curve $\mathcal{C}_F$ reads 
$$Bxy + Ax + Ay + A + B k + B + 1=0.$$
If $B=0$ then $A=1$, otherwise the curve contains an $\mathbb{F}_q$-rational line distinct from $x=y$.  In this case, $A=1$  would imply $C=0$ and $D=0$, a contradiction  since we assumed $\beta\neq 0$. 
If $B\neq 0$, then it is easily seen that if such a conic splits into two lines, then they are defined over $\mathbb{F}_q$, since its ideal points are both $\mathbb{F}_q$-rational and simple.

If $A+B+C+1=0$ then the polynomial $L(x,y)$ reads 
\begin{equation}\label{Conic}
B x^2 + Bxy  + By^2 + Ax+Ay + A + Bk + 1.
\end{equation}
Note that if $B=0$ then, as above, $A=1$, otherwise the line $Ax+Ay + A  + 1=0$ is $\mathbb{F}_q$-rational and different from $x=y$: again this would imply $C=D=0$, a contradiction. 

If $B\neq 0$ then $L(x,y)$ splits into two lines if and only if $A^2 + A B + B^2 k + B=0$. 
In this case the curve splits in 
$$\left(x+\omega y+\frac{A}{B\omega}\right)\left(x+\omega^2 y+\frac{A}{B}\omega\right)=0,$$
where $\omega^2+\omega+1=0$. These two lines are not defined over $\mathbb{F}_q$ if and only if $q$ is an odd power of $2$.

\section{The proof of the conjecture}\label{Sec:Proof}

In Proposition \ref{Proposizione:riassunto} we determined all the cases in which the curve $\mathcal{C}_F$ splits into components none of them is defined over $\mathbb{F}_q$. As already mentioned in  Introduction this is a necessary condition for $f_{\alpha,\beta}$ in \eqref{f_alfa_beta} to be a PP. 

In this section, we first prove that the conditions in Proposition \ref{Proposizione:riassunto} are equivalent to those of Conjecture \ref{Conjecture}. Recall that $k$ is  a fixed element in $\mathbb{F}_q$ with $Tr(k)=1$ and $i\in \mathbb{F}_{q^2}\setminus \mathbb{F}_q$ is such that $i^2=i+k$. Also, we denote by $\overline{k}$ the element such that $\overline{k}^4=k$. We consider $\alpha=A+iB$, $\beta=C+iD$, with $A,B,C,D \in \mathbb{F}_q$.

\begin{proposition}\label{prop1}
Case 1 in Proposition \ref{Proposizione:riassunto} happens if and only if 
\begin{equation}\label{CondCaso1}
\beta(1+\alpha^{q+1}+\beta^{q+1})+\alpha^{2q}=0, \quad \beta^{q+1}\neq 1, \qquad Tr\left( \frac{\beta^{q+1}}{\alpha^{q+1}}\right)=0, \quad \beta\in \mathbb{F}_q.
\end{equation}
\end{proposition}
\proof
The condition 
$\beta(1+\alpha^{q+1}+\beta^{q+1})+\alpha^{2q}=0$ 
is equivalent to
$$
\left\{
\begin{array}{l}
A^2C + A^2 + ABC + B^2Ck + B^2k + B^2 + C^3 + C^2D + CD^2k + C=0\\
D(B+D)^2k+A^2D + ABD + B^2 + C^2D + CD^2 + D=0\\
\end{array}
\right..
$$
Recall that we require $\alpha\beta\neq 0$. Therefore Conditions \eqref{CondCaso1} are equivalent to 
$$
\left\{
\begin{array}{l}
A^2C + A^2 + ABC + B^2Ck + B^2k + B^2 + C^3 + C^2D + CD^2k + C=0\\
D(B+D)^2k+A^2D + ABD + B^2 + C^2D + CD^2 + D=0\\
Tr\left( \frac{(C+iD)^{q+1}}{(A+iB)^{q+1}}\right)=0\\
(C+iD)^{q+1}\neq 1\\
D=0\\
C+iD\neq 0\\
A+iB\neq 0\\
\end{array}
\right.,
$$
that is

$$
\left\{
\begin{array}{l}
D=0\\
C^2\neq 1\\
B=0\\
Tr\left( \frac{C^2}{B^2}\right)=0\\
(A^2+C^2+C)(C+1)=0\\
A\neq 0\\
C\neq 0\\
\end{array}
\right. \iff
\left\{
\begin{array}{l}
D=0\\
C\neq 1\\
B=0\\
Tr\left( \frac{C}{B}\right)=0\\
A^2+C^2+C=0\\
A\neq 0\\
C\neq 0\\
\end{array}
\right. 
\iff
\left\{
\begin{array}{l}
A= \xi^2+\xi\\
B=0\\
C=\xi^2, \quad \xi \in \mathbb{F}_{q}\setminus \{0,1\}\\
D=0\\
Tr\left( \frac{\xi}{xi+1}\right)=0\\
\end{array}
\right..
$$

\endproof

\begin{proposition}\label{prop2}
Case 2 in Proposition \ref{Proposizione:riassunto} happens if and only if 
\begin{equation}\label{CondCaso2}
\beta(1+\alpha^{q+1}+\beta^{q+1})+\alpha^{2q}=0, \quad \beta^{q+1}\neq 1, \qquad Tr\left( \frac{\beta^{q+1}}{\alpha^{q+1}}\right)=0, \quad \beta\in \mathbb{F}_{q^2}\setminus\mathbb{F}_q.
\end{equation}

\end{proposition}
\proof
Conditions \eqref{CondCaso2} are equivalent to 
$$
\left\{
\begin{array}{l}
A^2C + A^2 + ABC + B^2Ck + B^2k + B^2 + C^3 + C^2D + CD^2k + C=0\\
D(B+D)^2k+A^2D + ABD + B^2 + C^2D + CD^2 + D=0\\
Tr\left( \frac{(C+iD)^{q+1}}{(A+iB)^{q+1}}\right)=0\\
C^2+CD+kD^2\neq 1\\
D\neq 0\\
C+iD\neq 0\\
A+iB\neq 0\\
\end{array}
\right..
$$
Note that $B=0$ would imply $A^2C + A^2 + C^3 + C^2D + CD^2k + C=0$ and $A^2D +  C^2D + CD^2 + D^3k + D=0$ from which we obtain $A^2D=0$, a contradiction. So $B\neq 0$.

Also, $B\neq D$. Suppose the contrary, then  $A^2C + A^2 + ABC + B^2k + B^2 + BC^2 + C^3 + C=0$ and $A^2B + AB^2 + B^2C + B^2 + BC^2 + B=B(A+C+1)(A+B+C+1)=0$. 
\begin{itemize}
\item If $C=A+1$ then $A^2 + AB + B^2k = B^2 + B$ and $C^2 + CD + D^2k=A^2+AB+B+B^2k+1=B^2+1$. Then 
$$0=Tr\left( \frac{\beta^{q+1}}{\alpha^{q+1}}\right)=Tr\left( \frac{C^2 + CD + D^2k}{A^2 + AB + B^2k}\right)=Tr\left( \frac{B^2+1}{B^2+B}\right)=Tr\left( 1+\frac{1}{B}\right).$$
On the other hand from $A^2 + AB + B^2k = B^2 + B$ we get $k=\frac{A^2}{B^2}+\frac{A}{B}+\frac{1}{B}+1$ and $1=Tr\left( k\right)=Tr\left(\frac{A^2}{B^2}+\frac{A}{B}+\frac{1}{B}+1 \right)=Tr\left(\frac{1}{B}+1 \right)$, a contradiction. 
\item If $C=A+B+1$ then $A^2 + AB + B^2k =  B$ and $C^2 + CD + D^2k=A^2+AB+B+B^2k+1=1$. $$0=Tr\left( \frac{\beta^{q+1}}{\alpha^{q+1}}\right)=Tr\left( \frac{C^2 + CD + D^2k}{A^2 + AB + B^2k}\right)=Tr\left( \frac{1}{B}\right)=Tr\left( \frac{1}{B}\right).$$
On the other hand from $A^2 + AB + B^2k = B$ we get $k=\frac{A^2}{B^2}+\frac{A}{B}+\frac{1}{B}$ and $1=Tr\left( k\right)=Tr\left(\frac{A^2}{B^2}+\frac{A}{B}+\frac{1}{B} \right)=Tr\left(\frac{1}{B} \right)$, a contradiction. 
\end{itemize}

This means that Conditions \eqref{CondCaso2} are equivalent to 
$$
\left\{
\begin{array}{l}
A^2C + A^2 + ABC + B^2Ck + B^2k + B^2 + C^3 + C^2D + CD^2k + C=0\\
D(B+D)^2k+A^2D + ABD + B^2 + C^2D + CD^2 + D=0\\
Tr\left( \frac{(C+iD)^{q+1}}{(A+iB)^{q+1}}\right)=0\\
C^2+CD+kD^2\neq 1\\
D\neq 0\\
B\neq 0\\
B\neq D\\
\end{array}
\right.,
$$
that is 
$$
\left\{
\begin{array}{l}
k=\frac{A^2D + ABD + B^2 + C^2D + CD^2 + D}{D(B+D)^2}\\
(AD + BC + BD + B)\cdot\\
\hspace{0.5 cm}(AD^2 + B^3 + BCD + BD^2 + BD)=0\\
Tr\left( \frac{(C+iD)^{q+1}}{(A+iB)^{q+1}}\right)=0\\
C^2+CD+kD^2\neq 1\\
D\neq 0\\
B\neq 0\\
B\neq D\\
\end{array}
\right..
$$

Note that $A=\frac{BC + BD + B}{D}$ would imply 
$$D(B+D)^2(C^2 + CD + D^2k + 1)=0$$
and this is a contradiction to $C^2 + CD + D^2k=\beta^{q+1}\neq 1$. So,

$$
\left\{
\begin{array}{l}
A=\frac{B^3 + BCD + BD^2 + BD}{D^2}\\
k=\frac{B^4 +C^2D^2 +CD^3 +D^2}{D^4}\\
Tr\left( \frac{C^2+CD+D^2k}{A^2+AB+B^2k}\right)=Tr\left( \frac{(B^2+D)^2/D^2}{B^2(B^2+D)/D^2}\right)=Tr\left( \frac{B^2+D}{B^2}\right)=0\\
C^2+CD+kD^2\neq 1\\
D\neq 0\\
B\neq 0\\
B\neq D\\
\end{array}
\right..
$$

Since 

$$1=Tr\left(k\right)=Tr\left(\frac{   B^4 + C^2D^2 + CD^3  + D^2 }{D^4}\right)=Tr\left(\frac{B}{D}+\frac{1}{D^2}\right),$$

the above conditions are equivalent to 

$$
\left\{
\begin{array}{l}
k=\frac{A^2D + ABD + B^2 + C^2D + CD^2 + D}{D(B+D)^2}\\
A=\frac{B^3 + BCD + BD^2 + BD}{D^2}\\
Tr\left( \frac{B^2+D}{B^2}+\frac{B}{D}+\frac{1}{D^2}\right)=1\\
C^2+CD+kD^2\neq 1\\
D\neq 0\\
B\neq 0\\
B\neq D\\
\end{array}
\right.,
$$
which are precisely Case 2 in Proposition \ref{Proposizione:riassunto}. 
\endproof

\begin{proposition}\label{prop3}
Case 3 in Proposition \ref{Proposizione:riassunto} happens if and only if
\begin{equation}\label{CondCaso3}
\beta=\alpha^{q-1}\in \mathbb{F}_q^{*}, \qquad Tr\left(1+\frac{1}{\alpha^{q+1}}\right)=0.
\end{equation}
 \end{proposition}
\proof
Conditions \eqref{CondCaso3} are equivalent to
$$
\left\{
\begin{array}{l}
D=0\\
(A+iB)C=A+(i+1)B\\
Tr\left(1+\frac{1}{A^2+AB+kB^2}\right)=0\\
C\neq 0\\
A+iB\neq 0\\
\end{array}
\right.
\iff
\left\{
\begin{array}{l}
D=0\\
A(C+1)+B=0\\
B(C+1)=0\\
Tr\left(1+\frac{1}{A^2+AB+kB^2}\right)=0\\
\end{array}
\right.
\iff
\left\{
\begin{array}{l}
A\neq 0\\
B=0\\
C=1\\
D=0\\
Tr\left(1+\frac{1}{A}\right)=0\\
\end{array}
\right..
$$
\endproof

\begin{proposition}\label{prop4}
Case 4 in Proposition \ref{Proposizione:riassunto} happens if and only if
\begin{equation}\label{CondCaso4}
\beta=\alpha^{q-1}\in \mathbb{F}_{q^2}\setminus \mathbb{F}_q, \qquad Tr\left(1+\frac{1}{\alpha^{q+1}}\right)=0, \qquad \alpha+\alpha^q\neq \beta+\beta^q.
\end{equation}
 \end{proposition}
\proof
Conditions \eqref{CondCaso4} are equivalent to
$$
\left\{
\begin{array}{l}
(A+iB)(C+iD)=A+B+iB\\
D\neq 0\\
Tr\left(1+\frac{1}{A^2+AB+kB^2}\right)=0\\
B\neq D\\
A+iB\neq 0\\
\end{array}
\right.
\iff
\left\{
\begin{array}{l}
AD+BC+BD+B=0\\
AC + A + BDk + B = 0\\
D\neq 0\\
Tr\left(1+\frac{1}{A^2+AB+kB^2}\right)=0\\
B\neq D\\
A+iB\neq 0\\
\end{array}
\right.$$
$$
\iff
\left\{
\begin{array}{l}
A=\frac{B(C+D+1)}{D}\\
B(C^2+CD+D^2k+1) = 0\\
D\neq 0\\
Tr\left(1+\frac{1}{A^2+AB+kB^2}\right)=0\\
B\neq D\\
A+iB\neq 0\\
\end{array}
\right..
$$
Clearly $B=0$ would imply $A=0$, a contradiction. So, the above conditions are equivalent to 
$$\left\{
\begin{array}{l}
A=\frac{B(C+D+1)}{D}\\
C^2+CD+D^2k+1 = 0\\
BD\neq 0\\
Tr\left(1+\frac{D}{B^2}\right)=0\\
B\neq D\\
\end{array}
\right..
$$
\endproof

\begin{proposition}\label{prop5}
Case 5 in Proposition \ref{Proposizione:riassunto} happens if and only if
\begin{equation}\label{CondCaso5}
\beta=\alpha^{q-1}\in \mathbb{F}_{q^2}\setminus \mathbb{F}_q, \qquad Tr\left(1+\frac{1}{\alpha^{q+1}}\right)=0, \qquad \alpha+\alpha^q= \beta+\beta^q.
\end{equation}
 \end{proposition}
\proof
Conditions \eqref{CondCaso4} are equivalent to
$$
\left\{
\begin{array}{l}
(A+iB)(C+iD)=A+B+iB\\
Tr\left(1+\frac{1}{A^2+AB+kB^2}\right)=0\\
B= D\neq 0\\
\end{array}
\right.
\iff
\left\{
\begin{array}{l}
AB+BC+B^2+B=0\\
AC + A + B^2k + B = 0\\
Tr\left(1+\frac{1}{A^2+AB+kB^2}\right)=0\\
B= D\neq0\\
\end{array}
\right.
$$
$$
\iff
\left\{
\begin{array}{l}
C=A+B+1\\
B= D\neq0\\
A^2+AB+B^2k+B= 0\\
Tr\left(1+\frac{1}{B}\right)=0\\
\end{array}
\right..
$$
Since $k=\frac{A^2+AB+B}{B^2}=\frac{A^2}{B^2}+\frac{A}{B}+\frac{1}{B}$ and then 
$$1=Tr(k)=Tr\left(\frac{A^2}{B^2}+\frac{A}{B}+\frac{1}{B}\right)=Tr\left(\frac{1}{B}\right)=Tr\left(1+\frac{1}{B}\right)+Tr(1)=Tr(1),$$
we have that the above conditions are equivalent to 
$$
\left\{
\begin{array}{l}
C=A+B+1\\
B= D\neq0\\
A^2+AB+B^2k+B= 0\\
q=2^{2s+1}\\
\end{array}
\right..
$$
\endproof

Now we are in position to prove Conjecture \ref{Conjecture}.

\begin{theorem}
Let $q=2^h$, $h\geq 3$. Let $f_{\alpha,\beta}(x)\in \mathbb{F}_{q^2}[x]$ as in \eqref{f_alfa_beta}. It permutes $\mathbb{F}_{q^2}$ if and only if  
\begin{enumerate}
\item $\beta=\alpha^{q-1}$ and $Tr\left(1+\frac{1}{\alpha^{q+1}}\right)=0$ or
\item $\beta(1+\alpha^{q+1}+\beta^{q+1})+\alpha^{2q}=0$, $\beta^{q+1}\neq 1$, and $Tr\left( \frac{\beta^{q+1}}{\alpha^{q+1}}\right)=0$.
\end{enumerate}
\end{theorem}
\proof
In Proposition \ref{Proposizione:riassunto} we completely described all the cases for which the curve $\mathcal{C}_F$ splits into absolutely irreducible components not defined over $\mathbb{F}_q$. By the Hasse-Weil Theorem if  the curve $\mathcal{C}_F$ has an absolutely irreducible component defined over $\mathbb{F}_q$ then it contains at least $q-2\sqrt{q}+1-4-4=q-2\sqrt{q}-7$ affine $\mathbb{F}_q$-rational points off the line $x=y$, since it has at most four ideal point and four points on $x=y$. So, for $q\geq 16$, if $f_{\alpha,\beta}$ is a PP of $\mathbb{F}_{q^2}$ then $\mathcal{C}_F$ has no affine $\mathbb{F}_q$-rational points off $x=y$ and therefore $\mathcal{C}_F$ splits completely into absolutely irreducible components not defined over $\mathbb{F}_q$, that is one of the conditions in Proposition \ref{Proposizione:riassunto} occurs. 

On the other hand, Propositions \ref{prop1}, \ref{prop2}, \ref{prop3}, \ref{prop4}, \ref{prop5} yield the equivalence between conditions in Proposition \ref{Proposizione:riassunto} and conditions in Conjecture \ref{Conjecture}. Since  \cite[Remark 2]{TZLH2017} points out that for $q=8$ Conjecture \ref{Conjecture} holds, the assertion follows.
\endproof

\section{Appendix}
In this section we present the MAGMA \cite{MAGMA} programs used in the previous sections. The computations involve the polynomial ring in the unknowns $x,y,A,B,C,D,i,k,a,b,c,d,e$ over the field $\mathbb{F}_2$. In some of our proofs we will need to investigate the solutions of a system of polynomial equations. In these situation an important tool we use is the resultant of two polynomials. 
We recall here some basic facts about the resultant. Given two polynomials  with coefficients in a field or in an integral domain $R$, say 
$$u(x) = a_mx^m + a_{m-1}x^{m-1} + \cdots + a_0, \qquad v(x) = b_nx^n + b_{n-1}x^{n-1} + \cdots + b_0,$$
with  $a_mb_n \neq 0$,  their resultant $Res(u,v)\in \mathbb{R}$ is the determinant of the  matrix
$$\left(
\begin{array}{cccccccccc}
a_m&a_{m-1}& \cdots & \cdots& \cdots&\cdots&a_0& 0& 0& 0\\
0&a_m&a_{m-1}& \cdots & \cdots& \cdots&\cdots&a_0&0 &0\\
&&\ddots&\ddots&&&&&\ddots\\
0&0 &0&a_m&a_{m-1}& \cdots & \cdots& \cdots&\cdots&a_0\\
b_n&b_{n-1}& \cdots & \cdots& b_0&0&0&0 &0 &0 \\
&\ddots&\ddots&&&\ddots\\
0&0&0& 0&0 &b_n&b_{n-1}& \cdots & \cdots& b_0\\
\end{array}
\right).$$
For a field $K$ and two polynomials
$F(x,y),\, G(x,y)\in K[x,y]$ of positive degree in $y$ we denote by $Res(F,G,y)$  the resultant of $F$ and $G$ with respect to $y$. It is the resultant of $F$ and $G$ when considered as polynomials in the single variable $y$ (that is, as elements in $R[y]$ with $R=K[x]$). In this case $Res(F,G,y)\in K[x]$ is in the ideal generated by $F$ and $G$, and therefore any pair $(a,b)$ with $F(a,b)=G(a,b)=0$ is such that $Res(F,G,y)(a)=0$;
see e.g. \cite[Prop 3.6.1]{COX}.

Recall that $k$ is an element of $\mathbb{F}_{2^m}$ of absolute trace $1$ and $i$ satisfies $i^2=i+k$. Also, $\alpha=A+iB$, $\beta=C+iD$, $\alpha^q=A+(i+1)B$, $\beta^q=C+(i+1)D$. Finally, the curve $\mathcal{C}_{F}$ has affine equation $(x+i+1)^3 (y+i+1)^3 h(x,y)/(x-y)=0$, where $h(x,y)$ equals
$$ \left(\alpha^q \left(\frac{x+i}{x+i+1}\right)^3+\left(\frac{x+i}{x+i+1}\right)^2+\beta^q\right) \left(\beta \left(\frac{y+i}{y+i+1}\right)^3+\left(\frac{y+i}{y+i+1}\right)+\alpha\right)$$
$$-\left(\alpha^q \left(\frac{y+i}{y+i+1}\right)^3+\left(\frac{y+i}{y+i+1}\right)^2+\beta^q\right) \left(\beta \left(\frac{x+i}{x+i+1}\right)^3+\left(\frac{x+i}{x+i+1}\right)+\alpha\right).$$

In what follows, the procedure ``FindCoefficients2" is used to quickly select all the coefficients of the monomials $x^iy^j$ in a bivariate polynomial ``pol". The procedure Substitution is used to substitute all the occurrences of a monomial $m$ in a polynomial ``pol" with the polynomial $p$.

{\footnotesize
\begin{verbatim}
FindCoefficients2 := function(pol,var1,var2)
	T := Terms(pol);
	Coeff := {};
	MAX1 := Degree(pol,var1);
	MAX2 := Degree(pol,var2);
	for i in [0..MAX1] do
		for j in [0..MAX2] do
			c := K!0;
			for t in T do
				if IsDivisibleBy(t,var1^i*var2^j) eq true and 
				   IsDivisibleBy(t,var1^i*var2^(j+1)) eq false  and 
				   IsDivisibleBy(t,var1^(i+1)*var2^j) eq false  then
					c := c+ K! (t/(var1^i*var2^j));
				end if;
			end for;
			if c ne 0 then
				Coeff := Coeff join {c};
				i,j,c;
			end if;
		end for;
	end for;
	return Coeff;
end function;


Substitution := function (pol, m, p)
	e := 0;
	New := K! pol;
	while e eq 0 do
		N := K!0;
		T := Terms(New);
		i:= 0;
		for t in T do
			if IsDivisibleBy(t,m) eq true then
				Q := K! (t/m);
				i := 1;
				N := K!(N + Q* p);
			else 
				N := K!(N + t);
			end if;
		end for;
		if i eq 0 then 
			return New;
		else	
			New := K!N;
		end if;	
	end while;
end function;
\end{verbatim}
}

{\footnotesize
\begin{verbatim}
K<x,y,A,B,C,D,i,k,a,b,c,d,e> := PolynomialRing(GF(2),13);

alpha := A+i*B;
alpha_q := A+(i+1)*B;
beta := C+i*D;
beta_q := C+(i+1)*D;
X := (x+i)/(x+i+1);
Y := (y+i)/(y+i+1);

Hxy := (alpha_q*X^3+X^2+beta_q)*(beta*Y^3+Y+alpha)
	-(alpha_q*Y^3+Y^2+beta_q)*(beta*X^3+X+alpha);

Curve := (x+i+1)^3*(y+i+1)^3*Hxy;
Curve  := Substitution(Curve,i^2,i+k);
\end{verbatim}
}

\subsection{$\gamma_{2,2}\neq 0$ and $\mathcal{C}_F$ splits into four lines}

Here we suppose that the curve splits into four lines. We want to find constrains on $A,B,C,D$. If $\mathcal{C}_F$ splits into four lines of type the polynomial 
$$Curve+(x+a) (x+b)(y+a)(y+b)(A^2 + A B + B^2 k + C^2 + C D + D^2k + D + 1)$$
must be the zero polynomial in $x$ and $y$. That it all its coefficients must vanish. We consider $B=D$ and $B\neq D$ separately. 
{\footnotesize
\begin{verbatim}
//CASE B=D 
PROD := (x+a)*(x+b)*(y+a)*(y+b);
CC := FindCoefficients2(Curve+
	(A^2 + A*B + B^2*k + C^2 + C*D + D^2*k + D + 1)*PROD,x,y);
CC := {Resultant(pol,B+D,D) : pol in CC};
{Factorization(pol) : pol in CC | pol ne 0};
//IN PARTICULAR a + b + 1=0 SINCE (A + C + 1)*(A + B + C + 1)!=0
p1 := a + b + 1;
CC1 := {Resultant(pol,p1,b) : pol in CC};
{Factorization(pol) : pol in CC1 | pol ne 0};
//IN PARTICULAR B=0
CC2 := {Resultant(pol,B,B) : pol in CC1};
{Factorization(pol) : pol in CC2 | pol ne 0};
//IN PARTICULAR A*k + A*a^2 + A*a + A + C*k + C*a^2
//+ C*a + C + k + a^2 + a=0 SINCE (A + C + 1)!=0
p3 := A*k + A*a^2 + A*a + A + C*k + C*a^2 + C*a + C + k + a^2 + a;
CC3 := {Resultant(pol,p3,a) : pol in CC2};
{Factorization(pol) : pol in CC3 | pol ne 0};
//THEN A^2 + C^2 + C=0 SINCE (A + C + 1)!=0
\end{verbatim}
}

{\footnotesize
\begin{verbatim}
//CASE B!=D
PROD := (x+a)*(x+b)*(y+a)*(y+b);
CC := FindCoefficients2(Curve+
	(A^2 + A*B + B^2*k + C^2 + C*D + D^2*k + D + 1)*PROD,x,y);
{Factorization(pol) : pol in CC | pol ne 0};
//IN PARTICULAR  p1 MUST VANISH
p1 := A^2*a + A^2*b + A^2 + A*B*a + A*B*b + A*B + B^2*k*a + B^2*k*b + B^2*k +
            B + C^2*a + C^2*b + C^2 + C*D*a + C*D*b + C*D + D^2*k*a + D^2*k*b + 
            D^2*k + D*a + D*b + a + b + 1;
CC1 := {Resultant(pol,p1,b) : pol in CC};
{Factorization(pol) : pol in CC1 | pol ne 0};
//IN PARTICULAR  p2 MUST VANISH
p2 := A^2*D + A*B*D + B^2*D*k + B^2 + C^2*D + C*D^2 + D^3*k + D;
CC2 := {Resultant(pol,p2,k) : pol in CC1};
{Factorization(pol) : pol in CC2 | pol ne 0};
//IN PARTICULAR  p3 MUST VANISH
p3 := A^2*D + A*B*D + A*D^2 + B^2*D*a^2 + B^2*D*a + B^2*D
	 + B^2 + B*D^2*a + C^2*D + D^3*a^2 + D^3 + D^2 + D;
CC3 := {Resultant(pol,p3,a) : pol in CC2};
{Factorization(pol) : pol in CC3 | pol ne 0};
//FINALLY (SINCE D!=0) A*D^2 + B^3 + B*C*D + B*D^2 + B*D=0
//HERE WE CHECK WHAT HAPPENS IF 
//A*D^2 + B^3 + B*C*D + B*D^2 + B*D=0 AND 
//A^2*D + A*B*D + B^2*D*k + B^2 + C^2*D + C*D^2 + D^3*k + D=0
Factorization(K!(D^9*Evaluate(Curve,
	[x,y,(B^3 + B*C*D + B*D^2 + B*D)/D^2,B,C,D,i,
	(B^4 + C^2*D^2 + C*D^3  + D^2)/D^4,a,b,c,d,e])));
\end{verbatim}
}

\subsection{$\gamma_{2,2}\neq 0$ and $\mathcal{C}_F$ splits into two conics switched by $(x,y) \mapsto(y,x)$}

Here we suppose that the curve splits into two absolutely irreducible conics.  If these two conics  are switched by $(x,y) \mapsto(y,x)$ then they have affine equation either
$$x^2+(a+i b)x+(a+(i+1)b)y+c=0, \qquad \textrm{ and } \qquad y^2+(a+(i+1)b) x+(a+i b) y+c=0,$$
for some $a,b,c \in \mathbb{F}_q$, 
or
$$xy+(a+i b)x+(a+(i+1)b)y+c=0, \qquad \textrm{ and } \qquad   xy+(a+(i+1)b) x+(a+i b) y+c=0,$$
for some $a,b,c \in \mathbb{F}_q$.

In the former case we immediately get a contradiction.
{\footnotesize
\begin{verbatim}
PROD := (x^2+(a+i*b)*x+(a+(i+1)*b)*y+c)*(y^2+(a+(i+1)*b)*x+(a+i*b)*y+c);
CC := FindCoefficients2(Curve+(A^2 + A*B + B^2*k + C^2 + C*D + D^2*k + D + 1)*PROD,x,y);
{Factorization(pol) : pol in CC | pol ne 0};
//THEN i*b + a + b=0, SINCE A^2 + A*B + B^2*k + C^2 + C*D + D^2*k + D + 1!=0
\end{verbatim}
}

The latter case is considered here.
{\footnotesize
\begin{verbatim}
PROD := (x*y+(a+i*b)*x+(a+(i+1)*b)*y+c)*(x*y+(a+(i+1)*b)*x+(a+i*b)*y+c);
PROD := Sostituzione(PROD,i^2,i+k,K);
CC := FindCoefficients2(Curve+(A^2 + A*B + B^2*k + C^2 + C*D + D^2*k + D + 1)*PROD,x,y);

//CASE B=D
{Factorization(pol) : pol in CC | pol ne 0};
CC := {Resultant(pol,B+D,D) : pol in CC};
{Factorization(pol) : pol in CC | pol ne 0};
//IN PARTICULAR p1 MUST VANISH SINCE (A + C + 1)*(A + B + C + 1)!=0
p1 := b + 1;
CC1 := {Resultant(pol,p1,b) : pol in CC};
{Factorization(pol) : pol in CC1 | pol ne 0};
//IN PARTICULAR B=0
CC2 := {Resultant(pol,B,B) : pol in CC1};
{Factorization(pol) : pol in CC2 | pol ne 0};
//IN PARTICULAR p3 MUST VANISH
p3 := A^2*k + A^2*c + A + C^2*k + C^2*c + k + c;
CC3 := {Resultant(pol,p3,c) : pol in CC2};
{Factorization(pol) : pol in CC3 | pol ne 0};
//IN PARTICULAR (C + 1)*(A + C + 1)*(A^2 + C^2 + C)=0
//IF (A^2 + C^2 + C) IT SPLITS INTO FOUR LINES
//THEN C=1 
CC4 := {Resultant(pol,C+1,C) : pol in CC3};
{Factorization(pol) : pol in CC4 | pol ne 0};
//FINALLY A*a^2 + A*a + A + 1=0

//CASE B!=D
PROD := (x*y+(a+i*b)*x+(a+(i+1)*b)*y+c)*(x*y+(a+(i+1)*b)*x+(a+i*b)*y+c);
PROD := Sostituzione(PROD,i^2,i+k,K);
CC := FindCoefficients2(Curve+(A^2 + A*B + B^2*k + C^2 + C*D + D^2*k + D + 1)*PROD,x,y);
{Factorization(pol) : pol in CC | pol ne 0};
//IN PARTICULAR p1=0
p1 := A^2*b + A^2 + A*B*b + A*B + B^2*k*b + B^2*k + B + C^2*b
	 + C^2 + C*D*b +C*D + D^2*k*b + D^2*k + D*b + b + 1;
CC1 := {Resultant(pol,p1,b) : pol in CC};
{Factorization(pol) : pol in CC1 | pol ne 0};
//IN PARTICULAR p2=0
p2 := A^2*k + A^2*c + A*B*k + A*B*c + A + B^2*k^2 + B^2*k*c + B*k + 
	B*c + C^2*k + C^2*c + C*D*k + C*D*c + D^2*k^2 + D^2*k*c + k + c;
CC2 := {Resultant(pol,p2,c) : pol in CC1};
{Factorization(pol) : pol in CC2 | pol ne 0};
//IN PARTICULAR p3=0
p3 := A^2*D + A*B*D + B^2*D*k + B^2 + C^2*D + C*D^2 + D^3*k + D;
CC3 := {Resultant(pol,p3,k) : pol in CC2};
{Factorization(pol) : pol in CC3 | pol ne 0};
//IN PARTICULAR 
//D(B + D)(A*D + B*C + B*D + B)(A*D^2 + B^3 + B*C*D + B*D^2 + B*D)=0
	//IF A*D + B*C + B*D + B=0
	CC4 := {Resultant(pol,A*D + B*C + B*D + B,A) : pol in CC3};
	{Factorization(pol) : pol in CC4 | pol ne 0};
	//THEN IN PARTICULAR 
	//B^2*a^2 + B^2*a + B^2 + B*C + B*D*a + B*D
	//        + B + C^2 + D^2*a^2 + D^2 + D+1=0

	//IF A*D^2 + B^3 + B*C*D + B*D^2 + B*D=0
	AA := (B^3 + B* C* D + B* D^2 + B* D)/D^2;
	kk := (B^4 + C^2 *D^2 + C *D^3  + D^2)/D^4;
	Factorization(K!(D^9*Evaluate(Curve,[x,y,AA,B,C,D,i,kk,a,b,c,d,e])));
\end{verbatim}
}

\subsection{$\gamma_{2,2}\neq 0$ and $\mathcal{C}_F$ splits into two conics both fixed by $(x,y) \mapsto(y,x)$}

Here we suppose that the curve splits into two absolutely irreducible conics. If these two conics are both fixed by $(x,y) \mapsto(y,x)$ then  they have affine equation
$$xy+(a+i b)x+(a+ib)y+c+id=0, \qquad \textrm{ and } \qquad xy+(a+(i+1) b)x+(a+(i+1)b)y+c+(i+1)d=0,$$
for some $a,b,c,d \in \mathbb{F}_q$.
{\footnotesize
\begin{verbatim}
PROD := (x*y+(a+i*b)*x+(a+i*b)*y+(c+i*d))*(x*y+(a+(i+1)*b)*x+(a+(i+1)*b)*y+(c+(i+1)*d));
PROD  := Substitution(PROD,i^2,i+k);
CC := FindCoefficients2(Curve+(A^2 + A*B + B^2*k + C^2 + C*D + D^2*k + D + 1)*PROD,x,y);
{Factorization(pol) : pol in CC | pol ne 0};
//IN PARTICULAR p1 MUST VANISH
p1 := A^2*d + A^2 + A*B*d + A*B + B^2*k*d + B^2*k 
	+ C^2*d + C^2 + C*D*d + C*D +D^2*k*d + D^2*k + D*d + d + 1;
CC1 := {Resultant(pol,p1,d) : pol in CC};
{Factorization(pol) : pol in CC1 | pol ne 0};
//IN PARTICULAR p2 MUST VANISH
p2 := A^2*b + A^2 + A*B*b + A*B + B^2*k*b + B^2*k + B + C^2*b
	 + C^2 + C*D*b + C*D + D^2*k*b + D^2*k + D*b + b + 1;
CC2 := {Resultant(pol,p2,b) : pol in CC1};
{Factorization(pol) : pol in CC2 | pol ne 0};
//IN PARTICULAR p3 MUST VANISH
p3 := A^2*k + A^2*a + A^2*c + A*B*k + A*B*a + A*B*c + A + B^2*k^2 + B^2*k*a + 
            B^2*k*c + B*k + B*c + C^2*k + C^2*a + C^2*c + C*D*k + C*D*a + C*D*c + 
            D^2*k^2 + D^2*k*a + D^2*k*c + k + a + c;
CC3 := {Resultant(pol,p3,c) : pol in CC2};
{Factorization(pol) : pol in CC3 | pol ne 0};
//IN PARTICULAR p4 MUST VANISH
p4 := A^4*a^2 + A^4*a + A^4 + A^3 + A^2*B^2*a^2 + A^2*B^2*a + A^2*B^2 + A^2*B*a
            + A^2*B + A^2*C + A^2*D*a + A^2 + A*B^2*k + A*B^2*a + A*B*C + A*B*D*a 
            + A*B + A*C^2 + A*C*D + A*D^2*k + A*D + A + B^4*k^2*a^2 + B^4*k^2*a + 
            B^4*k^2 + B^3*k*a + B^2*C*k + B^2*D*k*a + B*C^2*a + B*C*D*a + 
            B*D^2*k*a + B*D*a + B*a + C^4*a^2 + C^4*a + C^4 + C^3 + C^2*D^2*a^2 + 
            C^2*D^2*a + C^2*D^2 + C^2*D*a + C^2*D + C^2 + C*D^2*k + C*D^2*a + C + 
            D^4*k^2*a^2 + D^4*k^2*a + D^4*k^2 + D^3*k*a + D^2*a^2 + D^2 + D*a + D 
            + a^2 + a;
CC4 := {Resultant(pol,p4,a) : pol in CC3};
{Factorization(pol) : pol in CC4 | pol ne 0};
//IN PARTICULAR H_A_B_C_D MUST VANISH
H_A_B_C_D := A^6 + A^5*B + A^4*B^2*k + A^4*B^2 + A^4*C^2 + A^4*C*D + A^4*D^2*k + A^4
            + A^3*B^3 + A^2*B^4*k^2 + A^2*B^4*k + A^2*B^2*C^2 + A^2*B^2*C*D + 
            A^2*B^2*D^2*k + A^2*B^2 + A^2*C^4 + A^2*C^2*D^2 + A^2*C^2 + A^2*C*D 
            + A^2*D^4*k^2 + A^2*D^2*k + A^2*D + A*B^5*k^2 + A*B*C^4 + 
            A*B*C^2*D^2 + A*B*C^2 + A*B*C*D + A*B*D^4*k^2 + A*B*D^2*k + B^6*k^3 
            + B^4*C^2*k^2 + B^4*C*D*k^2 + B^4*D^2*k^3 + B^4*k^2 + B^2*C^4*k + 
            B^2*C^2*D^2*k + B^2*C^2*k + B^2*C*D*k + B^2*C + B^2*D^4*k^3 + 
            B^2*D^2*k^2 + B^2*D*k + B^2*D + C^6 + C^5*D + C^4*D^2*k + C^4*D^2 + 
            C^3*D^3 + C^2*D^4*k^2 + C^2*D^4*k + C^2 + C*D^5*k^2 + C*D + D^6*k^3 
            + D^2*k;
//THE FOLLOWING CONDITION IS EQUIVALENT TO  
Condition := Substitution(((a+b+i*b)^2+(c+d+i*d))*((a+i*b)^2+(c+i*d)),i^2,i+k);
//k^2*b^4 + k*b^2*d + k*d^2 + a^4 + a^2*b^2 + a^2*d + b^2*c + c^2 + c*d
//WE CHECK THAT Condition=0 if p1=p2=p3=p4=H_A_B_C_D=0
R1 := Resultant(Condition,p1,d);
R2 := Resultant(R1,p2,b);
R3 := Resultant(R2,p3,c);
R4 := Resultant(R3,p4,a);
Resultant(R4, H_A_B_C_D,k);
\end{verbatim}
}

\section*{Acknowledgments}

The author was partially supported by the Italian Ministero dell'Istruzione, dell'Universit\`a e della Ricerca (MIUR) and by the Gruppo Nazionale per le Strutture Algebriche, Geometriche e le loro Applicazioni (GNSAGA-INdAM).

\bibliographystyle{abbrv}

\begin{thebibliography}{10}

\bibitem{AGW2011} A. Akbary, D. Ghioca,  Q. Wang, On constructing permutations of finite fields, Finite Fields Appl. {\bf 17}, (2011) 51--67.

\bibitem{BG2017} D. Bartoli, M. Giulietti, Permutation polynomials, fractional polynomials, and algebraic curves, submitted.

\bibitem{MAGMA} W. Bosma, J. Cannon, and C. Playoust, The Magma algebra system. {I}. The user language, J. Symbolic Comput. {\bf 24}, (1997) 235--265. 

\bibitem{COX} D. Cox, J. Little, D. O'Shea, Ideals, Varieties, and Algorithms (Springer, 2007).

\bibitem{Dickson1896} L.E. Dickson, The analytic representation of substitutions on a power of a prime number of letters with a discussion of the linear group, Ann. Math. {\bf 11}, (1896) 65--120.

\bibitem{Hermite1863} Ch. Hermite, Sur les fonctions de sept lettres, C.R. Acad. Sci. Paris {\bf 57}, (1863) 750--757.

\bibitem{HirschBook} J.W.P. Hirschfeld, Projective geometries over finite fields, second edition,  Oxford Univ. Press, Oxford, (1998).

\bibitem{Hou2015} X. Hou, Permutation polynomials over finite fields--a survey of recent advances, Finite Fields Appl. {\bf 32}, (2015) 82--119.

\bibitem{MuPa} G.L. Mullen, and D. Panario, Handbook of Finite Fields (Chapman and Hall/CRC,
2013).

\bibitem{PL2001}  Y.H. Park, J.B. Lee, Permutation polynomials and group permutation polynomials, Bull.
Austral. Math. Soc. {\bf 63}, (2001) 67--74.

\bibitem{Sti}H. Stichtenoth: Algebraic Function Fields and Codes, 2nd Edition, Graduate Texts inMathematics, vol. 254, Springer, Berlin, 2009.

\bibitem{TZLH2017} Z. Tu, X. Zeng, C. Li, T. Helleseth, A Class of New Permutation Trinomials, Finite Fields and Their Applications, Finite Fields and Their Applications {\bf 50}, (2018) 178--195. 

\bibitem{YD2011}  P. Yuan, C. Ding, Permutation polynomials over finite fields from a powerful lemma,
Finite Fields Appl. {\bf 17}, (2011) 560 -- 574.

\bibitem{YD2014} P. Yuan and C. Ding, Further results on permutation polynomials over finite fields. Finite
Fields Appl. {\bf 27}, (2014) 88--103.

\bibitem{Zieve2009} M. Zieve. On some permutation polynomials over $\mathbb{F}_q$ of the form $x^rh(x^{(q-1)/d})$, Proc. Amer. Math. Soc. {\bf 137}, (2009) 2209--2216.



\end{thebibliography}
\end{document}